# BIFURCATIONS, SCHWARZIAN DERIVATIVE AND FEIGENBAUM CONSTANTS REVISITED

By Andrei Vieru
Paris, July 2006 - October 2007

**Abstract**
The main purpose of the this paper is to show that Feigenbaum δ constant is much more universal than believed. We'll speak about period-doubling processes in families of endomorphisms of the interval [0, 1] or of the semi-line [1, ∞[ depending on a parameter which is neither a term nor a factor; we'll first distinguish three cases:

$F_a(x) = [F(x)]^a$     one single (quadratic) minimum between 0 and 1
$F_a(x) = a^{g(x)}h(x)$     one single maximum (for *F*, *g* and *h*) between 0 and 1
$F_a(x) = [h(x)]^{a-g(x)}$     one single maximum (for *F*, *g* and *h*) between 0 and 1

Under reasonable conditions, in all these cases, the Feigenbaum δ constant appears. (For other cases and for the so-called *universality conjecture*, see Appendix 1.)
We formulate the so-called *parenthesis permeability hypothesis* – a conjecture that holds for all types of bifurcation (i.e. for *flip*, *fold*, *pitchfork* and *transcritical* bifurcations) – which states that under some conditions two or three different functions may have exactly the same bifurcation points (see also Appendix 2).
We propose a conjecture that considerably relaxes David Singer's conditions for endomorphism families to generate 'at most one stable orbit', showing that Feigenbaum δ constant appears also in some classes of functions that have more than one maximum and have positive Schwarzian in at least one sub-interval (see also Appendix 3).
Then, we consider some possible generalizations of the Feigenbaum δ constant. One of them was introduced by Feigenbaum himself, while studying the fixed points of the iterated function $f_\mu(x)=1 - \mu|x|^n$; we propose to generalize the Feigenbaum constant considering it as a derivative in the chaos point; then, one can study derivatives of higher order, partial and total derivatives, etc.
We propose at the end of the paper some unsolved or not yet studied problems related to generalized Feigenbaum constants.

# 1. The case $F_a(x) = [F(x)]^a$

The δ constant was thoroughly studied in relation with bifurcation processes generated by the increment of the parameter value the families of endomorphisms depend on. These parameters usually stand for a *term* or a *factor*. In some classical cases – including the celebrated logistic map

– the family of endomorphisms is supposed to satisfy:

$f_a(0) = f_a(1) = 0$      (1)
$f_a(x)$ – or $af(x)$) – has a single maximum between 0 and 1      (2)
The Schwarzian is supposed to be negative everywhere on [0, 1]      (3)

We started studying the behavior of endomorphism families where the parameter *a* is in an *exponent* position. The first family we looked at was $f_a(x)=x^{a/x}$ considered as a [1, ∞]→[1, ∞] map. One can easily see that for every $a > 0$

$f_a(1) = 1$      (4)
$f_a(x) \to 1$ when $x \to \infty$      (5)
$f_a(x)$ has a single maximum[1] in [1, ∞]      (6)

Considering different values of parameter *a*, we can find a period doubling process that leads to chaos, in many respects similar to the bifurcation process in the logistic map.

The $f_a(x)=x^{a/x}$ family displays some nice features:

- the first bifurcation point[2] is $e^2 = 7.3890560989…$ (see picture 1)

- the geometric means of the elements of every stable orbit equals the value of the parameter (probably a unique case)

- the 'Feigenbaum point' is situated very slightly lower[3] than $2e^2$

---

[1] namely for $x = e$
[2] Nevertheless, while writing these pages, we are unaware of the possibility to define *e* as the square root of the first bifurcation point of the family $f_a(x) = x^{a/x}$.
[3] Let's consider the functional Φ whose domain is the set of functions satisfying (4), (5) and (6), liable to lead to chaos through bifurcation process (by means of increasing the '*exponent*-parameter'), and whose range is the set of their Feigenbaum point ratios to their first bifurcation point. And let's consider a similar functional Λ, whose domain is the set of functions leading to chaos through increment of the '*factor*-parameter' and satisfying the above mentioned conditions (1) and (2). (We do not suppose (3) satisfied, because, as it will be shown later, this condition can be weakened.) 'Obviously', the ranges of these two functionals are bounded. But are their boundaries identical? Both ranges are 'obviously' included in ]1,2] but we conjecture that the upper boundary of the range of Φ is higher than the upper boundary of the range of Λ. In fact, we conjecture that the upper boundary of the range of Φ is reached precisely by the $f_a(x) = x^{a/x}$ family (roughly 14,77…/7.389…≈1.9989…)

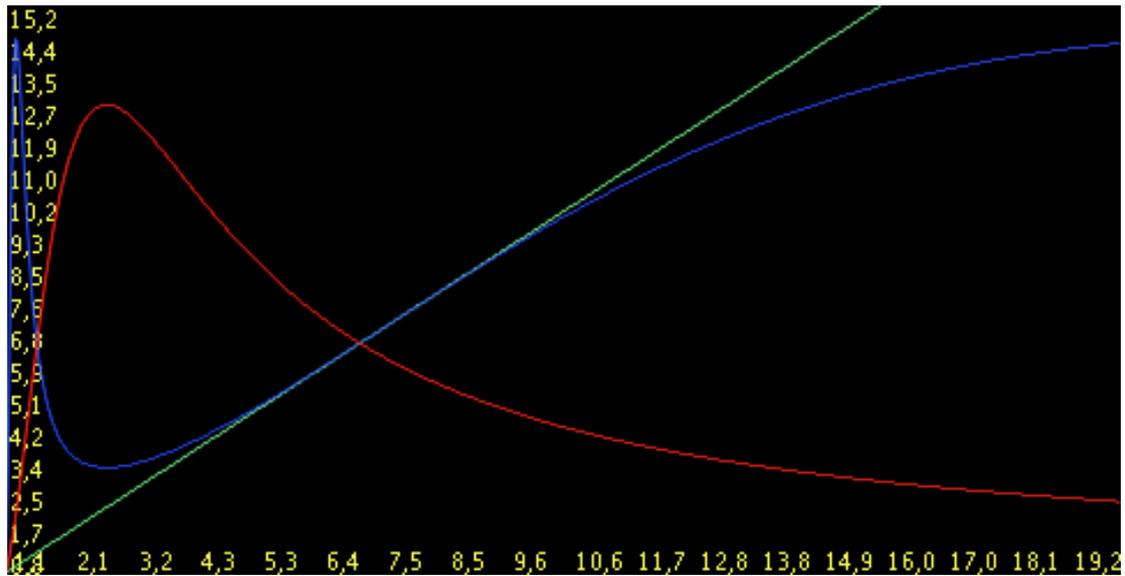

picture 1
$f_{7.389056}(x)=x^{\wedge}(e^2/x)$
$y=f^2_{7.389056}(x)$
and $y=x$ graph

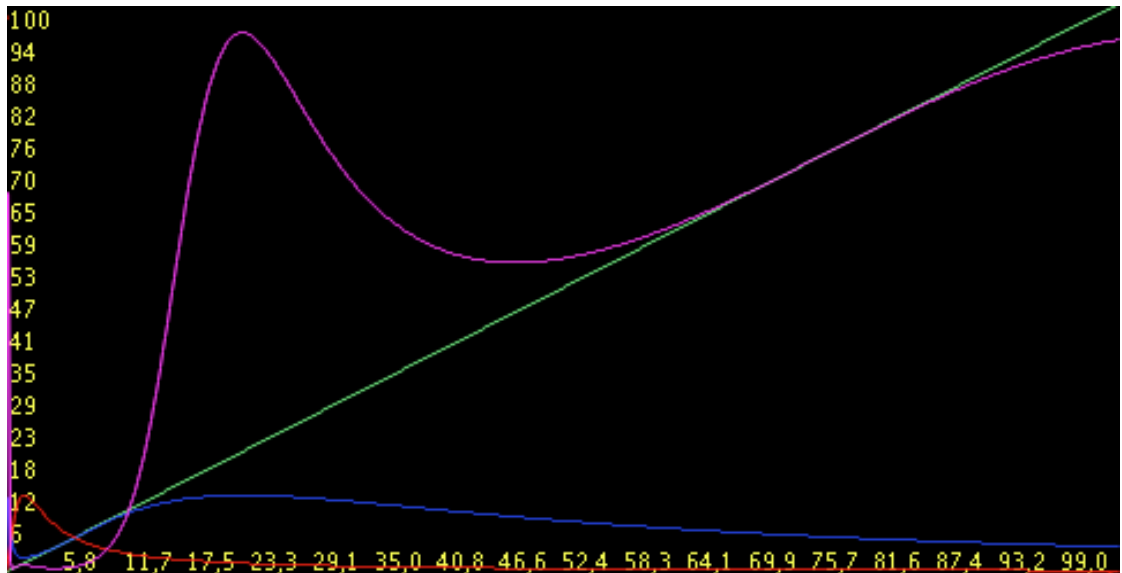

picture 2
$f_{7.389056}(x)=x^{\wedge}(e^2/x)$
$y=f^2_{7.389056}(x)$
$y=f^4_{12.5091}(x)$
and $y=x$ graph

Of course, writing $F(x) = [f(x^{-1})]^{-1}$ we may re-parameterize every map $f$ : $[1, \infty] \to [1, \infty]$ that satisfies conditions (4), (5) and (6), in order to obtain a map $F : [0, 1] \to [0, 1]$, such that

$$F(0) = F(1) = 1 \qquad (7)$$
$$F \text{ has one single minimum between 0 and 1} \qquad (8)$$

Obviously, $f_a$ and $F_a$ will have identical bifurcation points[4].

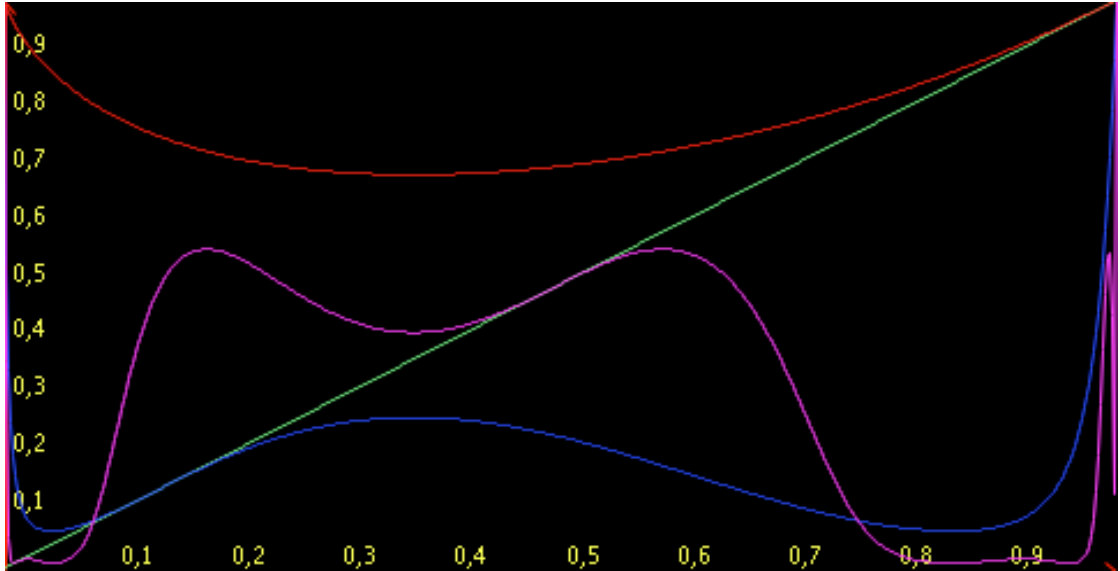

picture 3
$y=F_1(x) = x^{ax} = x^x$ (since $a = 1$)
$y=F^2_{7.389056}(x)$
$y=F^4_{12.5091}(x)$
$y=x$ graphs

Let us notice that the Schwarzian derivative of the function $f(x)=x^{1/x}$ is negative not everywhere on $[1, \infty]$. (In fact, it is positive[5] for every $x > 12.3944…$ )

Let us also notice that if the parameter stands for a *factor*, the Schwarzian – for obvious reasons – is the same for all members of a given family, while, when the parameter is an *exponent*, the Schwarzian – also for obvious reasons – is in general different for different functions of a given family. So, even if a lot of examples may bring it some credit, the hypothesis according to which $af(x) : [0, 1] \to [0, 1]$ leads to chaos through doubling period process if and only if $[1 – f(1 – x)]^a$ leads to chaos through doubling period process is highly unlikely, since it takes into account the principle of symmetry – in respect to the point (0.5, 0.5) – regardless of the behavior of the Schwarzian derivatives.

---

[4] This statement is true for any $f$ satisfying (4), (5) and (6) re-parameterized to $F$ satisfying (7) and (8), not only for our example. However, in our example – and probably only in it – the geometric means of any orbit generated by iterating $F_a$ equals $a^{-1}$.

[5] the Schwarzian of $f(x) = x^x$ is positive for every $x < 0.0806… = 1/12.3944…$ Moreover, $x^x$ has no finite derivative in 0 (assuming, using continuity, that $0^0 = 1$), but this fact does not impede on doubling period process leading to chaos.

Remarkably, calculations reveal that the δ constant appears in the example of the $f_a(x)=x^{a/x}$ family[6]. The Feigenbaum δ constant also appears in the case of $\xi_a(x)= [sin(\pi/x)+1]^a$ or in the case of $\psi_a(x)=[(x^2 + x – 1)/x^2]^a$, both considered as $[1, \infty]\to[1, \infty]$ maps[7]. So we believe, even if it is not *a priori* obvious, the Feigenbaum δ constant **does** appear everywhere families of endomorphisms satisfying (4), (5) and (6) – or (7) and (8) – lead to chaos through 'normal' bifurcation process. Moreover, this is true even if the Schwarzian is negative not everywhere. (We'll see further that Schwarzian may be 'negative not everywhere' in endomorphisms families of the interval [0, 1] that satisfy (1) and (2), and where the parameter they depend on stands for a *factor*.)

The role of Feigenbaum's α constant in the considered examples is not so clear. Examining maps on [0, 1], we could only see that in the case of $f_a(x) = x^{ax}$ we have a series of ramifications of the 'feigentree' through which passes the line[8] $y=e^{-1}$ instead of the line y = ½, as it happens with the logistic map[9].

---

[6] It is still difficult to assess the consequences of the fact that the *exponent*-parameter values may be as high as we wish. This was not the case of the logistic map (where $a \leq 4$) and other similar families of functions depending on a *factor*-parameter, where the set of its possible values is bounded. Anyway, we think the behavior of these families depending on an *exponent*-parameter beyond the chaos point deserves a glance: sometimes chaos slips into periodic orbits, then chaos is reached again through doubling period process, exactly as it happens in the case of the logistic map. For the $f_a(x) = x^{a/x}$ family, the geometric means of the elements of these stable orbits (beyond chaos point) continues to coincide with the parameter value. Does it make sense to consider that for chaos points this coincidence still holds?

[7] Let $\tau(k)$ be the ratio of the bifurcation point of rank $k$ in the above mentioned ($\psi_a$) family to the bifurcation point of rank $k$ in above mentioned $\xi_a$ family. Then, every pair of iterates $\psi_\alpha^k(x)$ and $\xi_{\tau(k)\alpha}^k(x)$ of order $k$, seem to 'quite well' approximate each other *on all* $[1, \infty]$. We have to confess that for us it is a quite disturbing phenomenon (even if, for big $k$ and α values, it might be less true)!

[8] we are unaware of the possibility to define $e^{-1}$ as the limit of the series of these elements of orbits just before bifurcation; these type of definitions may hide some circularity. Just before the *n*-th bifurcation point, the orbit has $2^{n-1}$ elements. If, in the orbits themselves, we consider the *littlest* element as the first, then the ranks, in the orbits, of the pairs of elements that serve to define the α constant are $2^{n-2}$ and $2^{n-1}$. Of course, in the case of the logistic map, if we consider the *biggest* element of each of these orbits as the first, then we'll find the elements used in order to define the α constant at places $2^{n-2}$ and $2^{n-1}$ respectively. But it is sometimes more convenient to

After a few bifurcations we found a ratio of approximately 2.559092… so we may still 'hope' to find the α constant. Curiously enough, calculating, for example, *all* ratios between all eight widths at 5th bifurcation point and all sixteen widths at 6th bifurcation point, we found the following pairs of ratios (6.393143…, 2.319988…) (2.721484.., 5.3667912) (6.333262…, 2.319112..) (2.442613…, 5.982122…) (6.327533…, 2.322344…) (2.559092…, 5.694513…) (5.524527…, 2.632882…) (2.336717…, 6.304994…).

We confess we do not know if the similarity of the four 'blue' numbers has or not any significance. We even don't know if the degree of similarity of their sums and products (smallest product 14.54…, biggest product 14.83…) is relevant or not, if this similarity will increase or not with every bifurcation point, if it will converge or not like other function families, etc.

## 2. THE PARENTHESIS PERMEABILITY CONJECTURE

Puzzling but true, two different[10] families of functions such as $\Psi_a(x)=\{\exp[-\sin(\pi x)]\}^a$ and $\Xi_a(x)=\exp[-\sin(\pi x^a)]$ have exactly the same bifurcation points. See below the first and the second approximate bifurcations points for both function families.

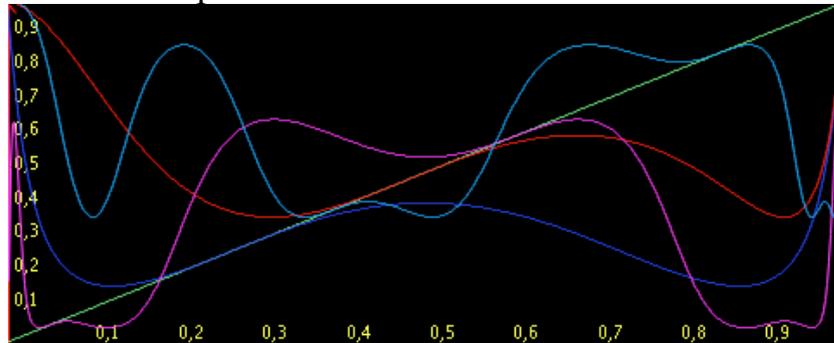

picture 4   $y=\Psi^2_{1.8}(x)$   $y=\Psi^4_{3.15}(x)$   $y=\Xi^2_{1.8}(x)$   $y=\Xi^4_{3.15}(x)$   $y=x$

---

consider the place of the elements in their orbits as we see them on the graph of the 'feigentree'. If, in the case of the logistic map, we prefer to consider the elements of a given orbit in the order of their magnitude (the biggest as the first, the littlest as the last), then, noting, after the *n*-th bifurcation, with A(*n*) the order of the smallest of the two elements considered when calculating the α constant, we'll have A(1)=1 and A(*n*+1)=2A(*n*) + (−1)$^{n+1}$. For maps with *exponent*-parameter satisfying (4) and (5), all happens in a similar way, only we'll consider the littlest element as the first, the biggest as the last, and we'll designate by A(*n*) the biggest of the two elements (of an orbit) that serve to calculate the α constant.

[9] We find the line $y = ½$, for example, in the case of $y = 1/(\sin(\pi x)+1)^a$

[10] The first is symmetric in respect to $x = 0.5$ vertical line, the second not.

Under reasonable conditions – maybe under no additional conditions at all – for every $c > 1$ and every $f : [0, 1] \rightarrow [0, 1]$ satisfying

$f(0) = f(1) = 0$           (1)
$f$ has a single maximum           (2)

$x \rightarrow c^{-af(x)}$ and $x \rightarrow c^{\wedge}-f(x^a)$ will always have the same bifurcation points.

We call this statement *parenthesis permeability hypothesis*. We also noticed that, if we replace $c$ by $x^r + s$ ($r > 0$, $s \geq 0$), the statement will not anymore be true, but 'big' differences of the $r$ values will induce 'very small' differences of the bifurcation points values: this is a kind of anti-chaotic behavior within a chaos related one. On the other hand, if $f$ satisfies (1) and (2), we'll have the same bifurcation points for $x \rightarrow [1 - f(x)]^a$ and for $x \rightarrow 1 - f(x^a)$.

This *parenthesis permeability* seems to be a very general phenomenon. It seems to take place even for 'bad' functions with 'bad' Schwarzians and 'bad' bifurcations, such as David Singer's example[11]: $f_a(x) = a(7.86x - 23.31x^2 + 28.75x^3 - 13.3x^4)$. The 'critical' values of the parameter in $x \rightarrow [1 - f(x)]^a$ and $x \rightarrow 1 - f(x^a)$ are the same, even if the orbits are not stable. The bifurcation values are also identical for $x \rightarrow af(x)$ and $x \rightarrow f(ax)$.

We can formulate an analogous hypothesis for some classes of maps where the parameter stands for a *factor*, namely for functions that have two critical points: one maximum and one minimum. If we take, for example, the function $\varphi(x) = x^3 - x$, we see that the families $x \rightarrow a\varphi(x)$ and $x \rightarrow \varphi(ax)$, although different, have exactly the same bifurcation points[12].

It seems that *parenthesis permeability conjecture* holds for all types of bifurcation (see also APPENDIX 2). Examining just the classical examples of the different types of bifurcation,

---

[11] *"Stable orbits and bifurcation of maps of the interval"* (SIAM vol. 35, No 2, September 1978)

[12] This $x \rightarrow \varphi(ax)$ family is particular: the interval whose endomorphism is considered shrinks when the value of the parameter increases. So it would be more correct to say 'family of endomorphisms defined on a family of intervals' (the family of intervals and the family of their endomorphisms depend, of course, on the same parameter). Chaos **is** reached, though on an interval that – as we just said – shrinks while the value of the parameter increases. Curiously, calculations reveal that the Feigenbaum δ constant still appears in this example.

we find, for *flip bifurcation*, that

$f_a(x) = a - x - x^2$ has exactly the same bifurcation points as $g_a(x) = -(x+a) - (x+a)^2$ ;

we find, for *fold bifurcation*, that

$f_a(x) = a - x^2$ has the same bifurcation points as $g_a(x) = -(a+x)^2$ ;

we also find, for *pitchfork bifurcation* that

$f_a(x) = ax - x^3$   $g_a(x) = a(x - x^3)$ and $h_a(x) = ax - (ax)^3$ have the same bifurcation points.

We find, for *transcritical bifurcation*, that

$f_a(x) = ax - x^2$, $g_a(x) = a(x - x^2)$ and $h_a(x) = ax - (ax)^2$ have the same bifurcation points.

Although this is not the subject of this paper, there is another interesting special case: the case of functions that have infinitely many critical points such as $\gamma_a(x) = a[\sin(2\pi x)/2 + x]$ considered as a **R**→**R** map. For $a < 1$, we have some interesting bifurcation values of parameter $a$ that still coincide with bifurcation values of the family $x \to \sin(2\pi ax)/2 + ax$. As stated below, some of the parameter values create *locally* stable orbits[13], along with *local* chaos. Doubling period of *locally* stable orbits involves here their complete displacement. For example, for $a \in (0.6049\ldots, 0.7901\ldots)$ we have an instable fixed point in 0 and two stable periodic orbits (with period two), symmetric in respect to 0, one attracting **R**- the other attracting **R**+. For $a \in (0.7901\ldots, 0.87201\ldots)$ we have an instable fixed point in 0 and two stable orbits with period 4 – one situated in [–1, 0], the other in [0, 1], one orbit symmetric to the other in respect to 0. For $a > 0.87201\ldots$ we have an instable fixed point in 0, two locally stable periodic orbits with period 4, one attracting [–1, 0] the other in [0, 1], where they are located, two locally stable orbits with period 8, one attracting [–2, –1] the other [1, 2] (where they are located), two locally stable periodic orbits with period 2, one attracting [–3, –2] the other in [2, 3] (where they are located), and two 'hyper-local' stable fixed points in 3.77459… and –3.77459…, both attracting some tiny neighborhoods of them…

For $a = 1$, **R** splits into unitary length intervals each of them attracted to a stable orbit with period 2 (on the verge to split into period 4). What happens when $a > 1$?

For example, if $a = 1.05$, we have a stable orbit with period 4, which attracts the interval ]0, ≈0.99…[. (The interval extremities are indicated roughly.)

But we also have a periodic orbit with period 10, which attracts [≈0.99…, ≈1.98…]. We also have chaos – no stable orbits at all – in the intervals [≈1.98…,

---

[13] These orbits attract, but not all points of the domain **R**.

≈2.96…] and [≈2.96…, ≈3.94…], an orbit with period 6 which attracts the interval [≈3.94…, ≈4.93…], chaos in the interval [≈4.93…, ≈5.91…], an orbit with period 4 which attracts the interval [≈5.91…, ≈6.89…] etc.

## 3. The case $F_a(x) = a^{g(x)}h(x)$

Let us also notice that a bifurcation process may take place when the parameter is neither an exponent, nor a mere factor. Let's consider the case when there is a factor near a function, but it depends on a parameter **and** on the function near which it is placed. Let's consider the case in which instead of having

$x \to af(x)$ we have $x \to g(a, f(x))f(x)$

For instance, it turns out that $x \to a^{\sin(\pi x)}\sin(\pi x)$ – viewed as a $[0, 1] \to [0, 1]$ mapping – really produces a bifurcation process leading to chaos. More generally, families of the type $f_a(x) = a^{g(x)}h(x)$ (where $g$ and $h$ both satisfy (1) and (2), and both have 'good' Schwarzian derivative[14]) lead to chaos through doubling period process. Here again, calculations reveal that we'll find the Feigenbaum δ constant.

## 4. The case $F_a(x) = [g(x)]^{a-h(x)}$

The Feigenbaum δ constant also appears in cases where it isn't the increment but the **decrement** of the parameter value that leads to chaos through bifurcation process.

Here is such an example, satisfying conditions (1) and (2):
$f_a(x) = [x(1-x)]^{a-x(1-x)}$

in which, once again, the parameter $a$ is neither a *factor* nor a 'pure' *exponent*[15]. The first bifurcation point, here, is roughly equal to 0.35…, the second bifurcation point is approximately 0.265…, then chaos is reached 'very quickly', *i.e.* before $a$ falls to 0.25, under which iteration is not anymore possible. More generally, the bifurcation points of $f_a(x) = g(x)^{a-h(x)}$ (where $g$ and $h$ both satisfy (1) and (2), reach their maximum – *necessarily smaller than* 1 – for the same $x$ value and both have 'good' Schwarzian derivative) seem to always make it appear the Feigenbaum δ constant.

---

[14] We'll explain below why we don't say, as everybody does, 'Schwarzian negative everywhere'.

[15] Function families like $f_a(x) = [x(1-x)]^a$ do not bifurcate

# 5. SUFFICIENT CONDITIONS FOR AN ENDOMORPHISM FAMILY TO GENERATE CHAOS THROUGH DOUBLING-PERIOD PROCESSES

Considering the classical case when the parameter stands for a factor, let us notice that no one of the well-known David Singer's conditions for a function to bifurcate is necessary. David Singer qualifies these conditions as rather 'restrictive', so we would like to reformulate them in order to relax them as much as possible.
Let us, first of all, rewrite one of Singer's important results in this field:
THEOREM (David Singer[16])
Let $G$ be the set of smooth endomorphisms $F$ of $[0, 1]$ satisfying
  (a) $F(0) = F(1) = 0$
  (b) $F$ has a unique critical point $c$ in $[0, 1]$
  (c) $F$'s Schwarzian negative everywhere
 Then for any $F$ in $G$ there is at most one stable orbit in $(0, 1)$ ; if it exists it is the ω limit set of $c$. »

First of all, if we really assume condition (a) and condition (b) then, we state that the condition of 'Schwarzian negative everywhere' is not required; we only need a Schwarzian negative everywhere in some critical point's neighborhood, which changes it's sign *at most once* in $]c, 1]$ and *at most once* in $[0, c[$.
We shall not prove this statement, but we'll give two different examples, that might help someone 'feel' or 'see' (pictures 5, 6, 7) why it is so[17]:

(I)   $x \rightarrow a[1 - x^x(1 - x)^{1-x}]$
(which does bifurcate, leads do Feigenbaum δ constant[18], but has positive Schwarzian in $[0, 0.081...[$ and $]0.918..., 1]$ intervals and, moreover, has no finite derivative at all in 0 and in 1; the Schwarzian tends to +∞ as $x \rightarrow 1$ and as $x \rightarrow 0$)

---

[16] *"Stable orbits and bifurcation of maps of the interval"* (SIAM vol. 35, No 2, September 1978)

[17] Let us stress it again: according to my hypothesis, the whole trouble with the Schwarzian is not to see it changing its sign as we go toward the interval extremities, but to see it changing it's sign *again* before reaching 0 or 1.

[18] Another example of a function family without finite derivative in 0 and 1: $f(x) = [x(1 - x)]^{1/2}$. It 'normally' bifurcates, leads to Feigenbaum constant, so the condition of 'unique critical point', along with 'thrice derivability' are too strong.

(II)  $x \to a(-1.55 x^4 + 4.34x^3 - 4.56 x^2 + 1.77x)$

(which has positive Schwarzian in the interval ]0.7, 1], finite derivative in 0 and 1. Calculations show that this family also engenders the Feigenbaum δ constant.)

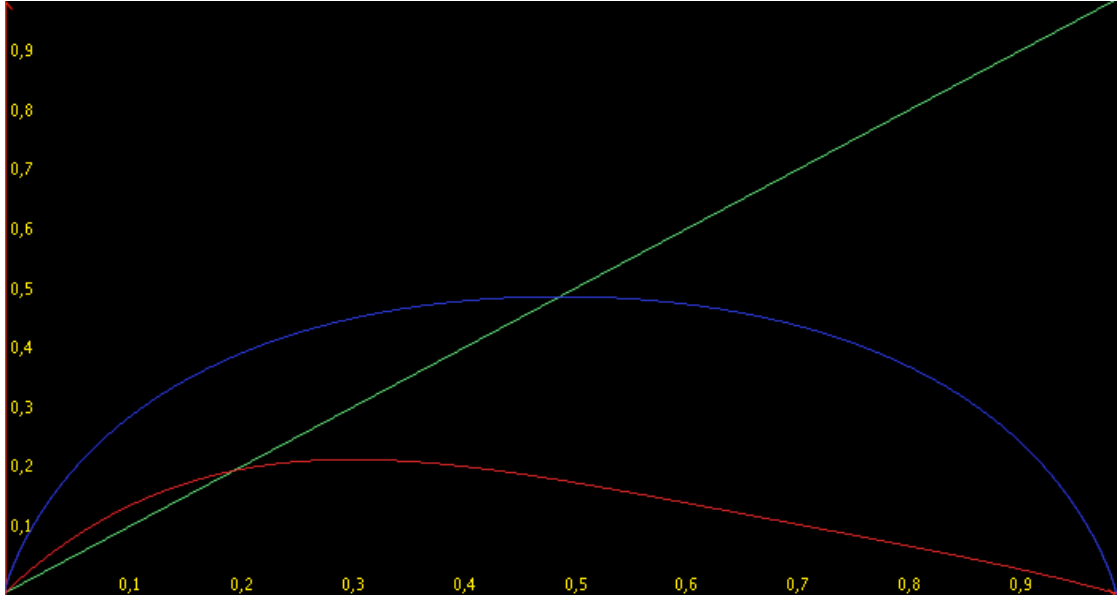

picture 5
$f(x) = 1 - x^x(1-x)^{1-x}$  $g(x) = -1.55 x^4 + 4.34x^3 - 4.56 x^2 + 1.77x$  and
$y = x$ graphs

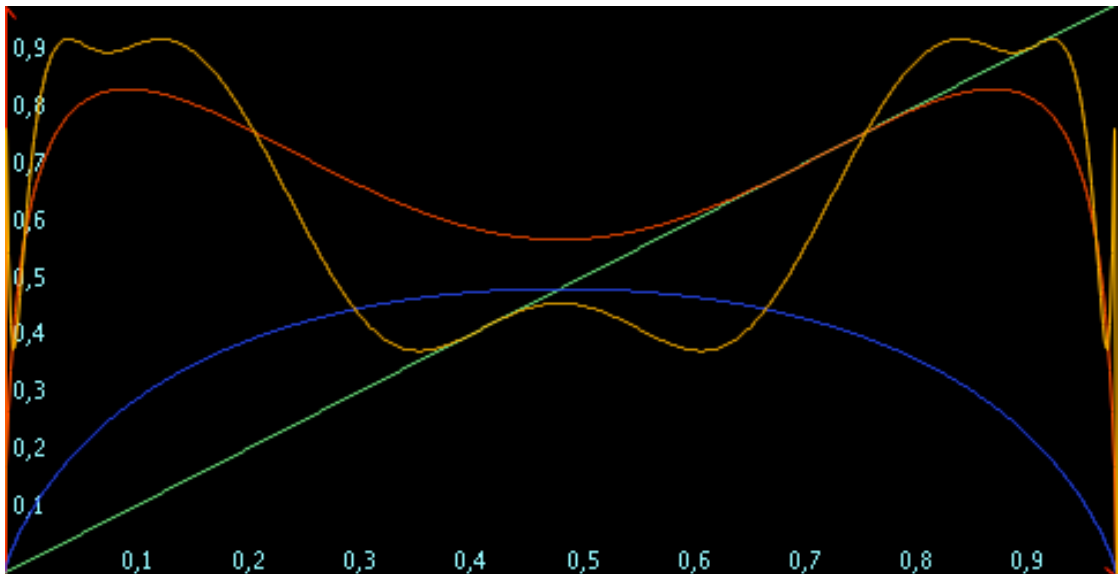

picture 6
$f(x) = 1 - x^x(1-x)^{1-x}$
$y = f^2_{1.7}(x)$
$y = f^4_{1.875}(x)$ and
$y = x$ graphs

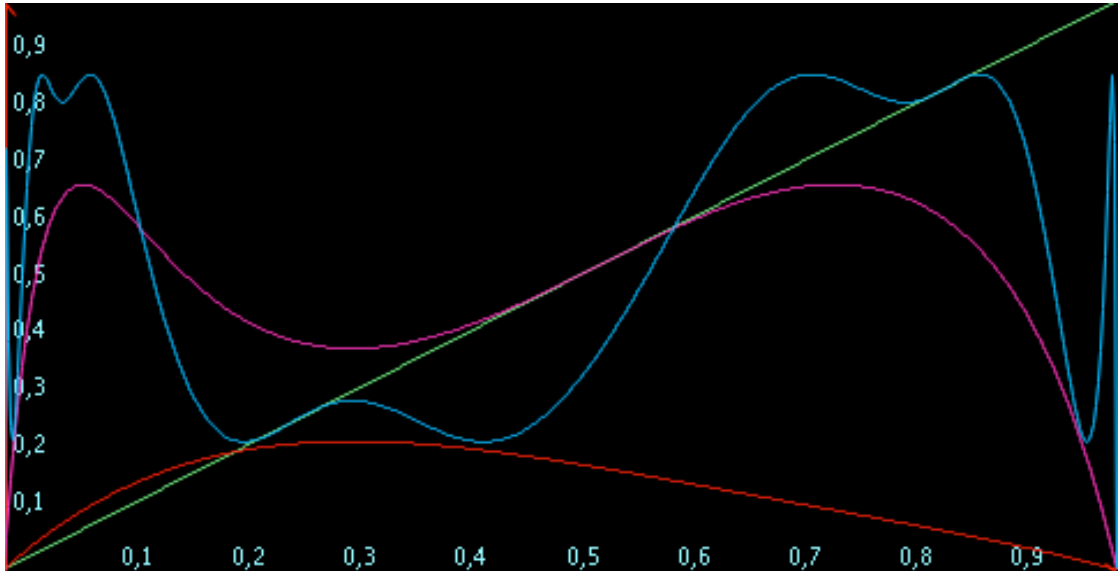

picture 7
$g(x) = -1.55 x^4 + 4.34 x^3 - 4.56 x^2 + 1.77x$
$y = g^2_{3.001}(x)$
$y = g^4_{3.86}(x)$ and
$y = x$ graphs

Now, if we consider the second condition, we immediately see that bifurcations may occur when we have several local maxima.
$x \to a(-x^8 + 4x^3 - 5x^2 + 2x)$ has two unequal local maxima, negative Schwarzian everywhere, and does bifurcate (see picture 8). Here again, the Feigenbaum δ constant appears!

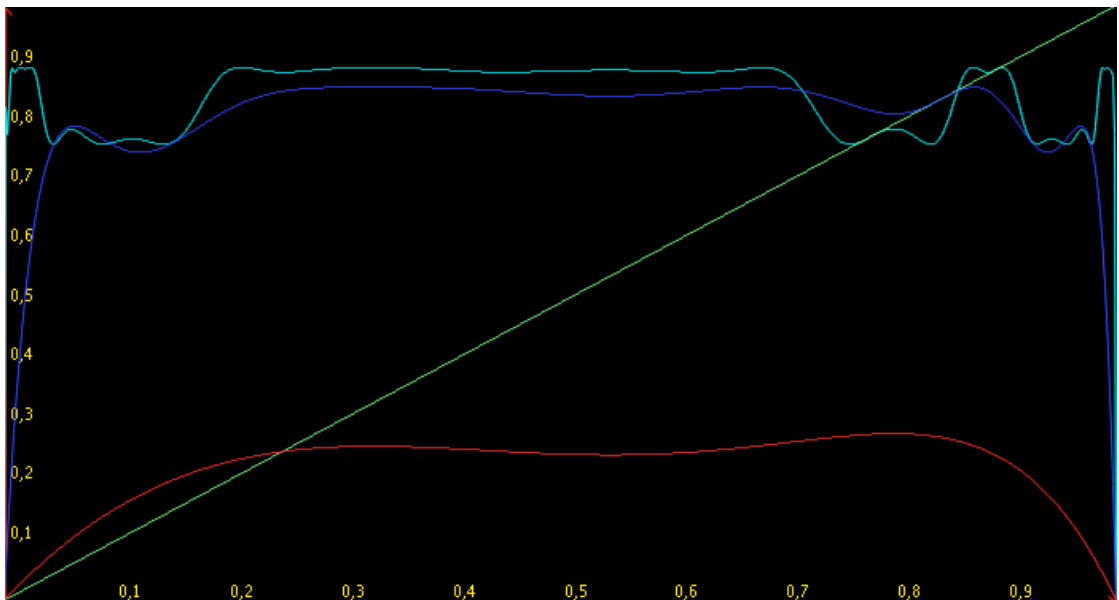

picture 8
$F(x) = -x^8 + 4x^3 - 5x^2 + 2x$   $y = F^2_{3.0781}(x)$   $y = F^4_{3.19746}(x)$  and $y=x$ graphs

Here are three examples of functions – each of them have three local maxima – that lead to chaos following the doubling period process.

$$x \to 9x(1-x)[1-3x(1-x)] + x^4(1-x^4) + 0.25 \times [-(2x-1)^2 + 1]^{44}$$

$$x \to 2.75 \times (-x^8 + 4x^3 - 5x^2 + 2x) + 0.05 \times [-(2x-1)^2 + 1]^{42}$$

$$x \to 2.75 \times (-x^8 + 4x^3 - 5x^2 + 2x) + 0.03 \times [-(2x-1)^2 + 1]^{42}$$

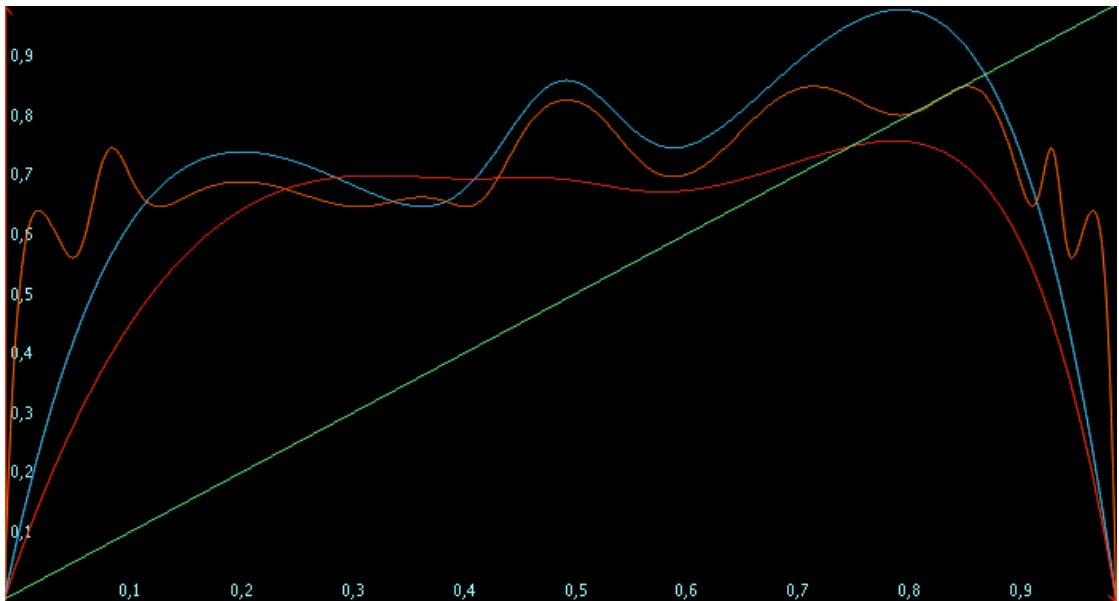

picture 9
$f(x) = 2.75 \times (-x^8 + 4x^3 - 5x^2 + 2x) + 0.03 \times [-(2x-1)^2 + 1]^{42}$
$g(x) = 9x(1-x)[1-3x(1-x)] + x^4(1-x^4) + 0.25 \times [-(2x-1)^2 + 1]^{44}$
$y = g^2_{0.87}(x)$ and $y = x$ graphs

## 6. THE SEVERAL MAXIMA BIFURCATION CONJECTURE:

Let $f_a(x) = af(x)$ be a family of endomorphisms of the interval $[0, 1]$ and let $\{x_1, x_2, \ldots, x_n\}$ be the set of $x$ values in which local quadratic maxima are attained. Let us assume that $x_1 < x_2 < \ldots < x_n$

If:
a) $\forall i < n \ f(x_i) \leq f(x_n)$
b) the Schwarzian is negative everywhere in $[x_1, x_n]$
c) the Schwarzian changes its sign at most once in $]x_n, 1]$

then the family engenders stable orbits and chaos is reached through period doubling process. Here again the Feigenbaum constant always appears (see also APPENDIX 3)!

Please note that when we have more than one maximum, the Schwarzian may change its sign in the open interval $]0, x_1[$ *more than once!*

The following function has two local maxima and Schwarzian positive on an interval strictly included in $]0, x_1[$:
$x \to a \times \{0.15 \times [1 - (-2x - 1)^4] + 0.5 \times \{1 - [2(1 - x)^8 - 1]^4\} + 2.4 \times x^2(1 - x^2)\}$

Nevertheless, it seems to reach chaos through bifurcation process (see picture 10). Here again, the Feigenbaum δ constant appears!

To be complete, let's notice that even the condition (1) is unnecessary to reach chaos through bifurcations. However, we need to have $f(0) > f(1)$.

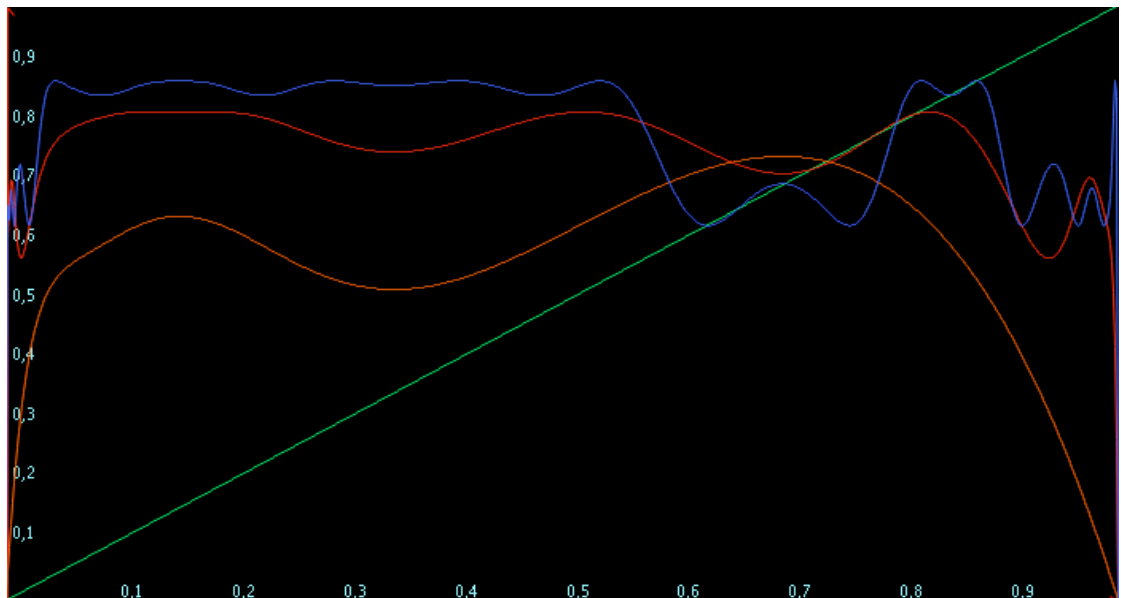

picture 10
$f(x) = \{0.15 \times [1 - (-2x - 1)^4] + 0.5 \times \{1 - [2(1 - x)^8 - 1]^4\} + 2.4 \times x^2(1 - x^2)\}$,
$y = f^2_{1.1}(x)$,
$y = f^4_{1.17}(x)$   and
$y = x$ graphs

# 7. ON GENERALIZED FEIGENBAUM CONSTANTS

The question of the generalization of the Feigenbaum δ constant may be raised in at least three different directions.

### I

One of them has been already studied by Mitchell Feigenbaum himself, namely the case when the maximum is not necessarily quadratic. On this issue, we have to ask:

1) Did anybody noticed that even[19] 'Feigenvalues' are situated on some smooth curve, while odd[20] 'Feigenvalues' (arising in the study of functions with a maximum of odd degree, *i.e.* non-analytical) are situated on *another* smooth curve?

So do both the row of differences between consecutive even and odd Feigenvalues and the row of differences between consecutive odd and even Feigenvalues (both are asymptotes to O$x$ axis).

We have no knowledge of these formulas – if they exist – of all these curves, but we would like to stress that the search of one single formula for both even and odd Feigenvalues leads only to wrong and ugly results.

2) Since non-analytical maxima are studied, why only of integer degree[21], and why always equal to the right and to the left – we mean why not of a 'hybrid' degree (of different degrees leftward and rightward[22])?

3) What happens if we have two critical points between 0 and 1, one quadratic, the other of degree 4 (or vice-versa)? We only can provide an example for further examination. For $a$=1.19 the 2$^{nd}$ iterate already shows a split not into a period 4 (?!) stable orbit, but into 2 sub-intervals;

---

[19] Feigenvalues that arise when in the iterated function $f_\mu(x)=1-\mu|x|^n$ $n$ is even

[20] Feigenvalues that arise when in the iterated function $f_\mu(x)=1-\mu|x|^n$ $n$ is odd

[21] in fact some non integer degree where computed when calculating the feigenvalue limit when in the iterated functions $f_\mu(x)=1-\mu|x|^n$ $n \to 1$.

[22] For example what is the feigenvalue of the iterated function $f_\mu(x)=1-\mu|x|^3$ if $x<0$ and $f_\mu(x)=1-\mu|x|^8$ if $x \geq 0$ ?

each of them is attracted by a period 2 orbit (see picture 11). After the following bifurcations points – the first four points of bifurcation are, roughly, 1.18336, 1.357, 1.3871, 1. 39154 – one can see that the part of the graph related to the degree 4 maximum takes over the whole: either in terms of Feigenbaum constant (of degree 4, i.e. 7.2846862171...) or in terms of location of the orbits (see also APPENDIX 3).

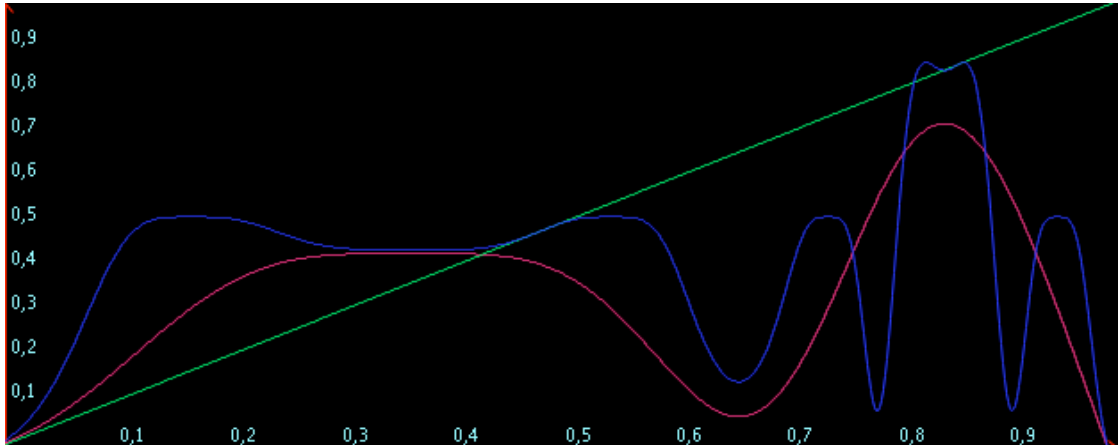

picture 11
$f_a(x)=a((0.9(\sin(0.9\pi(\sin(\pi(2(0.76x + 0.42) – 2)) + 2(0.76x + 0.42) –2)) + 0.9(\sin(\pi(2(0.42 + 0.76x) – 2)) + 2(0.76x + 0.42) – 2)) + 2 + 3(1 –16(0.76x + 0.42 – 0.7)^4) –1.5)/3 –0.6 )$   ($a=1$)
$y=f^2_{1.19}(x)$ and
$y=x$ graphs

II

Let $\{b_1, b_2, …\}$ be the set of bifurcation points of a given family. If we consider the function $b_n \to b_{n-1}$ (considered on $\{b_2, b_3, …\}$), then the Feigenbaum δ constant might be defined as the derivative of this function in cluster point $b_\infty$. So, when we have several parameters, one can ask what would be the 'partial' or 'total' derivative in every cluster point situated on the cluster curve. We already can state that the Feigenbaum constant appears sometimes as a 'partial derivative' i.e. changing the value of only one parameter, while keeping the others fixed[23].

One can study the simplest case of two-parameters map family
1°)   $af(x) + bg(x)$
(We suppose $f$ and $g$ both satisfy at least (1). Some interesting problems may arise in the study of this very simple form: what if $f$ and $g$ both have some quadratic critical points, whereas the whole form has one

---

[23] see also our paper 'Generalized Iteration, Catastrophes and Generalized Sharkovsky's Ordering' **arXiv:0801.3755 [math.DS]**

single maximum of degree 4. What if $f$ and $g$ both have 'bad' Schwarzians but the entire form has a 'good' one?)

For example (see picture 12), let us consider the function $y = a \times 1.5625 \times (0.25 - (x - 0.5)^2) + b \times (0.25 - (2.5 \times (0.25 - (x - 0.5)^2) - 0.5)^2)$

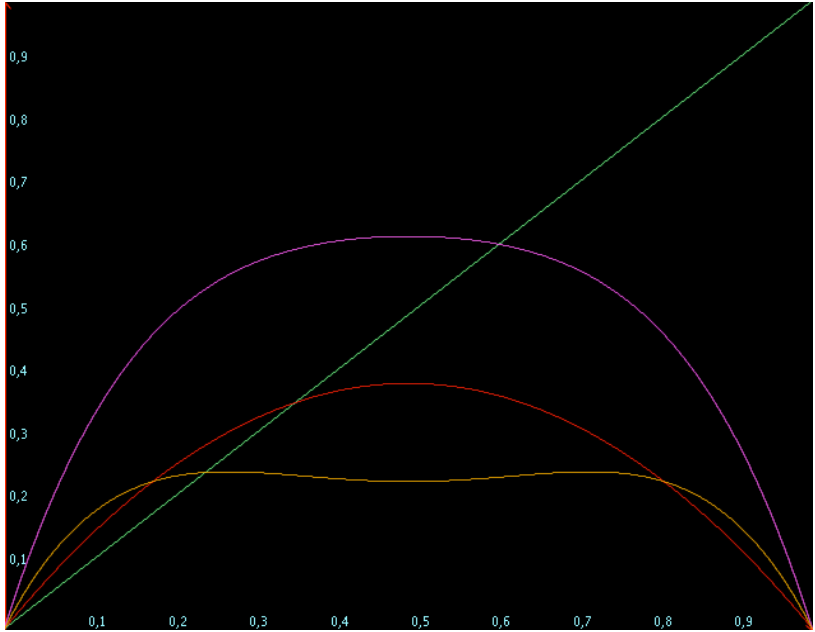

picture 12
$f(x)=1.5625\times(0.25 - (x - 0.5)^2)$ (one maximum of degree 2 in 0.5)
$g(x)=(0.25-(2.5\times(0.25-(x-0.5)^2)-0.5)^2)$ (three critical points one of which is a minimum of degree 2 in 0.5)
$h(x)= f(x) + g(x) = 1.5625\times(0.25 - (x - 0.5)^2) + (0.25-(2.5\times(0.25-(x-0.5)^2)-0.5)^2)$ (one maximum of degree 4)
$y=x$

Let's consider the bifurcation points of the function family depending on two parameters $h_{a,b}(x)= af(x) + bg(x)$. If we keep $a=1$ then the rate at which the distance between critical $b$ values[24] decreases will converge to some Feigenvalue (possibly the Feigenbaum constant for degree 4 maxima). If we keep $b=1$ then the rate at which the distance between critical $a$ values[25] decreases will converge to some other Feigenvalue (possibly 4.669… , the Feigenvalue for quadratic maxima). The reader is invited to check…

---

[24] i.e. for which $h_{a,b}(x)$ bifurcates
[25] for which $h_{a,b}(x)$ bifurcates

Here are the graphs (pictures 13 and 14) for these iterated functions. They seem suggestive.

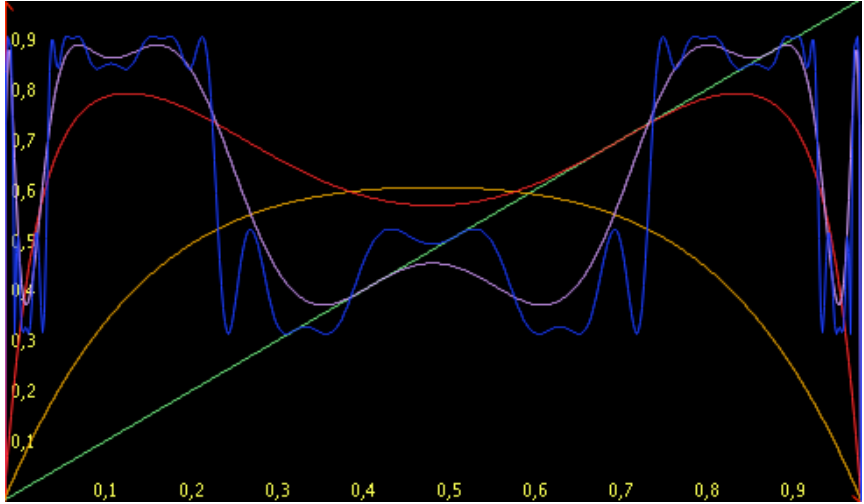

picture 13
$h(x) = f(x) + g(x) = 1.5625 \times (0.25 - (x - 0.5)^2) + (0.25 - (2.5 \times (0.25 - (x - 0.5)^2) - 0.5)^2)$
$h^2_{1.48,\ 1}(x) = 1.48f(1.48f(x) + g(x)) + g(1.48f(x) + g(x))$
$h^4_{1.725,\ 1}(x)$
$h^8_{1.77,\ 1}(x)$
$y = x$

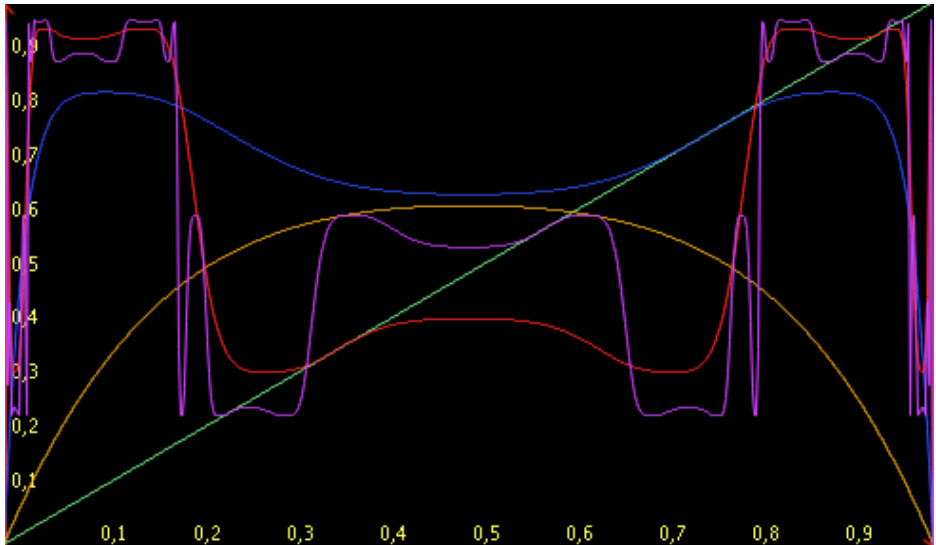

picture 14
$h(x) = f(x) + g(x)$
$h^2_{1,\ 1.89}(x) = f(f(x) + 1.89g(x)) + 1.89g(f(x) + 1.89g(x))$
$h^4_{1,\ 2.392}(x)$    $h^8_{1,\ 2.455}(x)$    $y = x$

Although in these pictures the bifurcations points are calculated very roughly, their ratios limits seem very differently oriented.
$(1.725 - 1.48)/(1.77 - 1.725) = 5.44444\ldots$ while
$(2.392 - 1.89)/(2.455 - 2.392) = 7.9682539\ldots$

Besides, picture 13 looks typically like the beginning of a bifurcation process of a function family with a quadratic maximum, while picture 14 looks typically like the beginning of a bifurcation process of a function family with a degree 4 maximum...

But what can be said about all other directions, all other straight lines of the plan $(a, b)$ on which we chose $a$ and $b$ values? What if neither $a$ nor $b$ are fixed? What about the rate at which, following all other directions, the distance between critical $(a, b)$ points will decrease? Obviously, when $a=b$ we'll find the feigenvalue for degree 4 $(7.28468\ldots)$, since $h_{a,b}(x) = af(x) + bg(x)$ has a maximum of degree 4.

Some analogous questions arise for some other map families of the following forms:

2°)   $af(x)bg(x)$ (at least $f$ satisfying (1) and (2))
3°)   $af(x)[g(x)]^b$ ($f$ satisfying (1) and (2), while $g$ satisfies (7) and (8))
4°)   $af(x)bg(x)$ ($f$ satisfying (1) and (2), while $g$ satisfies (7) and (8), a special case of 2°))

### III

If we consider the sequence of values that converges to the value of the derivative in $b_\infty$ (I mean the sequence formed by the terms $c_n = (b_\infty - b_{n-1})/(b_\infty - b_n)$ converging to Feigenbaum $\delta$ constant), we can define, as in II, an analogous function $c_n \to c_{n-1}$ and ask what would be the derivative in cluster point $c_\infty = \delta = 4.669\ldots$

If this derivative exists – let us call it, inappropriately, second derivative – we might ask what is it like, and of course try to see if the next *derivatives* do exist or not.

When are they all equal to the Feigenbaum constant itself? Always? Never? Sometimes?

## 8. WHEN THE FEIGENBAUM CONSTANT DOES NOT APPEAR

We already have given the example of $\gamma_a(x) = a[\sin(2\pi x)/2 + x]$ considered as a R→R map, where we could not find any $\delta$ constant.

In more usual contexts let's notice that for any family $f_a(x)$ of endomorphisms of the interval that reaches chaos and make it appear the

Feigenbaum δ constant, a 're-parameterized' family of endomorphisms of the interval of the form $\theta(a)f(x)$ (or of the form $f(x)^{\theta(a)}$) may not engender the Feigenbaum constant if, for the value $\lambda$ of the parameter $a$ for which chaos is reached, we have $\theta'(\lambda)=0$ or $\theta'(\lambda)=\infty$

Andrei Vieru
andreivieru2004@yahoo.fr

**Acknowledgements**
We express our deep gratitude Vlad Vieru and Dmitry Zotov for computer programming and to Robert Vinograd, Terente Robert and Sergiu Klainerman for useful discussions.

## APPENDIX 1 (the 'universality conjecture')

We have already considered the cases:
$F_a(x) = [F(x)]^a$    one single minimum between 0 and 1
$F_a(x) = a^{g(x)}h(x)$    one single maximum (for $F$, $g$ and $h$) between 0 and 1
$F_a(x) = [h(x)]^{a-g(x)}$    one single maximum (for $F$, $g$ and $h$) between 0 and 1

Of course there are other cases. Let us only consider three pairs of twin cases. They may be instructive[26]:

1°) $F_a(x) = -\text{Log}(af(x))/[\text{Log}(g(x))]^2$
2°) $F_a(x) = -\text{Log}(f(x))/[\text{Log}(ag(x))]^2 =$
$-\text{Log}(f(x))/\{[\text{Log}(a)]^2 + 2\text{Log}(a)\text{Log}(g(x)) + [\text{Log}(g(x))]^2\}$
($f(0)=f(1)=g(0)=g(1)=0$, one single quadratic maximum smaller than 1 for $f$, $g$; one single quadratic maximum for $F$ in both cases, for sufficiently low – respectively high – parameter $a$ values in 1°) – respectively in 2°))

In the first case – under the assumption[27] that the *decrement* of the parameter $a$ value generates chaos through bifurcation process – calculations quickly reveal the presence of the Feigenbaum $\delta$ constant.
Although in the second case the parameter $a$ has a different status, calculations also quickly reveal the presence of the Feigenbaum $\delta$ constant (under the assumption[26] that the *increment* of the parameter $a$ value generates chaos through bifurcation process).

3°) $F_a(x) = -\text{Log}(af(x))/[\text{Log}(g(x))]^2$
4°) $F_a(x) = -\text{Log}(f(x))/[\text{Log}(ag(x))]^2$
($f(0)=f(1)=1$, one single quadratic minimum; $g(0)=g(1)=0$ one single quadratic maximum; $F_a(0)=F_a(1)=0$ (using continuity) one single quadratic maximum)

---

[26] and not only because the exponent in the denominator may be either any even strictly positive integer or any odd integer ≥ 3 – but then, we have to drop the minus sign in the numerator

[27] this assumption may be a superfluous caution: when $f$ and/or $g$ have 'bad' Schwarzian it may happen that for different $a$ values $F_a$ has two different maxima in [0, 1]. Moreover it may happen that for some $a$ values, we may have $x_1<x_2$ and $F_a(x_1)<F_a(x_2)$, while for some other $a$ values we may have $x_3<x_4$ and $F_a(x_3)>F_a(x_4)$ ($x_1$, $x_2$, $x_3$ and $x_4$ designate values of the variable $x$ for which maxima are reached). Example: $g(x)$ is of the form $r_1(7.86x - 23.31x^2 + 28.75x^3 - 13.3x^4)$ while $f(x)$ is of the form $r_2(1-x^x)$. However, once we assume the uniqueness of the maximum for all $a$ values for which $F_a(x)$ remains an endomorphism, we could not find examples of families of functions of types 1°), 2°), 3°) or 4°) for which chaos and Feigenbaum constant do not appear.

In 3°), calculations quickly reveal the presence of the Feigenbaum δ constant (under the assumption that the *decrement* of the parameter *a* value generates chaos through bifurcation process).

In 4°), calculations quickly reveal the presence of the Feigenbaum δ constant (under the assumption that the *increment* of the parameter *a* value generates chaos through bifurcation process).

**THE UNIVERSALITY CONJECTURE:**

For any family of endomorphisms of the interval [0, 1] of any form $y = F_a(x) = F(a, x)$ such as
(1) $F_a(0)=F_a(1)=0$ for all *a* and such as
(2) there is a single quadratic maximum in [0, 1],
there are only two possibilities: either changes in the parameter *a* values do not engender chaos, or changes in the parameter *a* values do engender chaos *and then, the limit of the ratio between two consecutive differences in the sequences of bifurcation points equals δ=4.669…*

## APPENDIX 2 (the parenthesis permeability conjecture)

$F(x)+a$ and $F(x+a)$ have always the same bifurcations points (their set may be empty)

$aF(x)$ and $F(ax)$ have always the same bifurcations points (their set may be empty)

$[F(x)]^a$ and $F(x^a)$ have always the same bifurcations points (their set may be empty)

## APPENDIX 3 (about the several maxima conjecture)

Let's consider a two maxima example:
$$y= aF(x)=a\{1.2x^{7.9}(1-x^{7.9})+(1-x)^2[1-(1-x)^2]\}$$

here is its graph together with its forth iteration (for *a*=3.03):

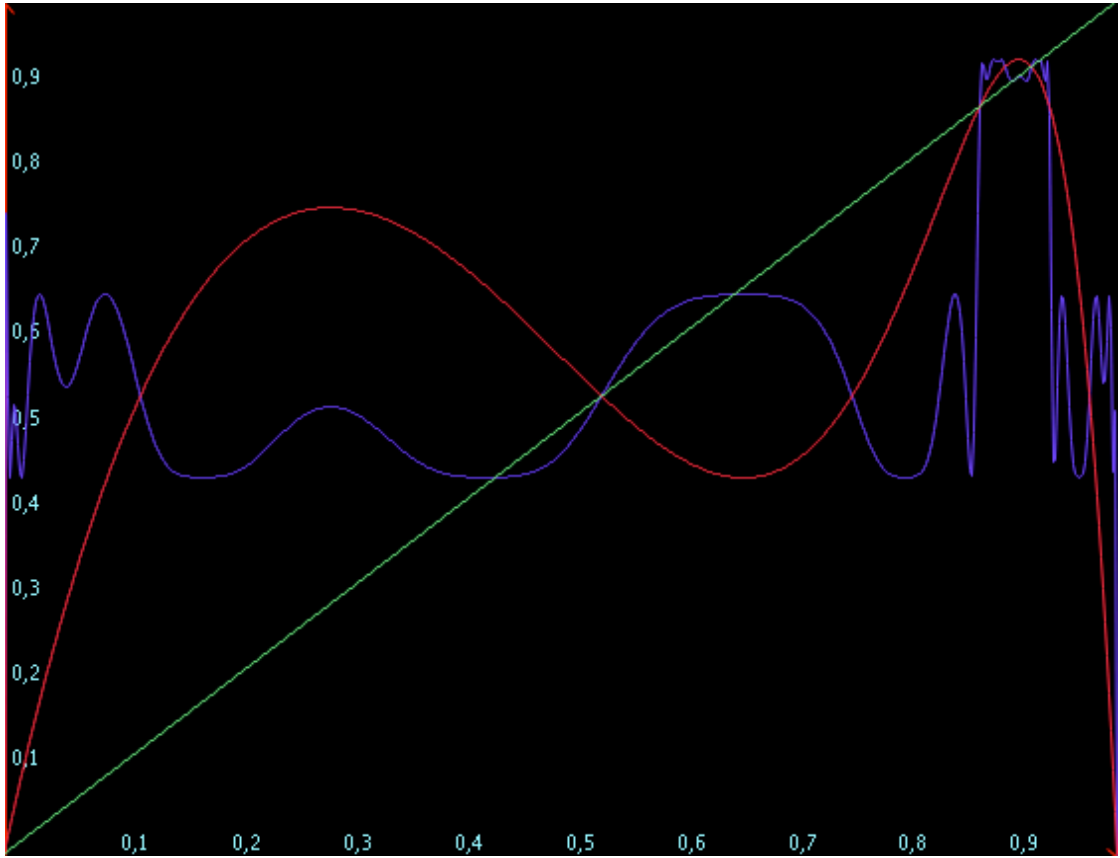

as one can see, the blue graph shows that the iteration of this endomorphism 'splits' the interval [0, 1] into subintervals: chaos is reached but only on the subinterval [≈0.87, ≈0.945]. Things happen like if the right maximum bifurcates faster than the left one (so, in some right part of the picture, local chaos will be reached before the left part will reach the second bifurcation).

In picture 11, we already have seen an example in which the left maximum bifurcates faster. There is no surprise in this phenomenon: in picture 11, the left maximum is of degree 4 (so it bifurcates at a rate tending to 7.2846862171...), while the right maximum is of degree 2 (so it bifurcates at a rate tending to 4.669201…).

We propose the reader to answer the following question: under which conditions a two quadratic maxima endomorphism engenders stable orbits that attract the whole interval [0, 1], as in the examples showed in picture 8 and in picture 10? Let $x_1$ and $x_2$ be the two points in which maxima are attained. Are there to positive numbers $\lambda$ and $\mu$, such as if $|x_1-x_2|\leq\delta$ and/or $f_a(x_2)-f_a(x_1)\leq\mu$, then stable orbits encountered in the bifurcation process that leads to chaos attract the whole interval [0, 1]? We believe thing are more complicated: the value of the minimum attained between the two maxima plays an important role too.